\documentclass[plainhead, 12pt,twoside]{article}

\usepackage{amssymb}
\makeatletter
\renewcommand\section{\@startsection{section}{1}{\z@}%
                                  {-3.5ex \@plus -1ex \@minus -.2ex}%
                                  {2.3ex \@plus.2ex}%
                                  {\normalfont\bfseries}}
\makeatother

\renewcommand{\leq}{\leqslant}
\renewcommand{\geq}{\geqslant}

\newcommand{\Rkg}{R_{KG}}
\newcommand{\Rkc}{R_{KC}}

\newcommand{\Rpnu}{R_{p^\nu}}
\newcommand{\Rpm}{R_{p^m}}

\newcommand{\Rpq}{R_{pq}}

\begin{document}
\parindent=0.3in
\parskip=0in
\baselineskip=18pt plus 1pt

\vspace*{0.1in} \noindent {\Large \bf Adams operations on the Green
ring of a cyclic group of prime-power order}$^\bigstar$

\vspace{0.2in} \noindent
{\large R. M. Bryant, Marianne Johnson}$^\ast$\footnote
{\hspace{-0.2in}$^\bigstar$ Work supported by EPSRC
Standard Research Grant EP/G024898/1.\\
\hspace*{0.08in}$^\ast\,$ Corresponding author.\\
\hspace*{0.2in}{\it E-mail addresses}: roger.bryant@manchester.ac.uk
(R. M. Bryant),\\
\hspace*{0.2in}marianne.johnson@maths.manchester.ac.uk (M. Johnson).}

\vspace{0.2in} \noindent
{\it School of Mathematics, University of Manchester,
Manchester M13 9PL, UK}

\vspace{0.4in} \noindent {\bf Abstract:} We consider the
Green ring $R_{KC}$ for a cyclic $p$-group $C$ over a field $K$ of
prime characteristic $p$ and determine the Adams operations $\psi^n$ in
the case where $n$ is not divisible by $p$. This gives information on
the decomposition into indecomposables of exterior powers and
symmetric powers of $KC$-modules.

\vspace{0.2in} \noindent {\bf Keywords:} Adams
operation, cyclic $p$-group, exterior power,
symmetric power

\vspace{0.4in} \baselineskip=20pt plus 1pt \noindent {\bf 1.
Introduction}
\renewcommand{\theequation}{1.\arabic{equation}}
\setcounter{equation}{0}

\noindent Let $C$ be a cyclic group of order $p^\nu$,
where $p$ is a prime and $\nu \geq 1$, and let $K$ be a field of
characteristic $p$. It is
well known that there are, up to isomorphism, exactly $p^\nu$
indecomposable $KC$-modules, and these can be written as $V_1, V_2,
\ldots, V_{p^\nu}$, where $V_r$ has dimension $r$, for $r=1, \ldots,
p^\nu$. The exterior powers $\Lambda^n(V_r)$ and symmetric powers
$S^n(V_r)$ have been studied intermittently
for more than thirty years. Some of the
main contributions have been by Almkvist and Fossum~\cite{AF},
Kouwenhoven~\cite{Ko}, Hughes and Kemper~\cite{HK}, Gow and
Laffey~\cite{GL} and Symonds~\cite{Sy}. The main aim has been to describe
$\Lambda^n(V_r)$ and $S^n(V_r)$, up to isomorphism, as direct sums
of indecomposable modules. An explicit formula is probably not
feasible, but one can look for a recursive description, so that, for
example, $\Lambda^n(V_r)$ is described in terms of exterior powers
$\Lambda^m(V_j)$ where $m<n$ or $j<r$.
The case $\nu=1$ was settled in~\cite{AF}, although further
information was provided by a number of people in subsequent
papers. However, for $\nu>1$, the problem remains open in general.

It is helpful to work in the Green ring (or representation ring)
$\Rkc$. This consists of all formal $\mathbb{Z}$-linear combinations
of $V_1, V_2, \ldots, V_{p^\nu}$, with addition defined in the
obvious way and multiplication coming from the decomposition of
tensor products into indecomposables. Finite-dimensional
$KC$-modules may be regarded, up to isomorphism, as elements of
$\Rkc$. This ring was first studied in detail by Green~\cite{Gr} in
1962, and he gave recursive formulae that implicitly describe
multiplication in $\Rkc$. Improved formulae and algorithms were
subsequently given by several other people: see, for
example,~\cite{McF,Ra,Re,Sr}.

In this paper we study the Adams operations $\psi_\Lambda^n$ and
$\psi_S^n$, for $n \geq 1$, following the treatment of these
in~\cite{Br}. Both $\psi_\Lambda^n$ and $\psi_S^n$ are
$\mathbb{Z}$-linear maps from $\Rkc$ to $\Rkc$. Furthermore,
$\Lambda^n(V_r)$ is given in $\mathbb{Q} \otimes_{\mathbb{Z}} \Rkc$
as a polynomial in $\psi_\Lambda^1(V_r), \ldots,
\psi_\Lambda^n(V_r)$. For example,
\begin{equation}
\Lambda^2(V_r) = \textstyle\frac{1}{2}(\psi_\Lambda^1(V_r)^2
- \psi_\Lambda^2(V_r)),
\end{equation}
where $\psi_\Lambda^1(V_r) = V_r$. Similarly, $S^n(V_r)$ is given
as a polynomial in $\psi_S^1(V_r), \ldots, \psi_S^n(V_r)$.

The main results of this paper determine $\psi_\Lambda^n(V_r)$ and
$\psi_S^n(V_r)$ for $n$ not divisible by $p$. Thus our results could
be used to determine $\Lambda^n(V_r)$ and $S^n(V_r)$ for $n<p$. For
$n$ not divisible by $p$, it is known (see \cite{Br}) that
$\psi_\Lambda^n=\psi_S^n$. Thus, in this case, we write $\psi^n$,
where $\psi^n=\psi_\Lambda^n=\psi_S^n$. In Section 3 we establish
the periodicity of these Adams operations (namely, $\psi^n=\psi^{n+2p}$)
and a
symmetry property (namely, $\psi^n=\psi^{2p-n}$ for $n=1, \ldots,
p-1$). We also prove a result (Proposition 3.6) that
generalises the ``reciprocity theorem'' of Gow and
Laffey~\cite[Theorem 1]{GL}. Most of the results of Section 3 extend
work for $\nu = 1$ by
Almkvist~\cite{Al} and Kouwenhoven~\cite{Ko}.

Our first main result (Theorem 4.7) describes $\psi^n(V_r)$
recursively in terms of the values $\psi^n(V_j)$ for $j<r$. This is
a simple recursion that enables $\psi^n(V_r)$ to be calculated in a
straightforward way by elementary arithmetic, and (strangely enough)
the recursion does not require any ability to multiply within
$\Rkc$.

One can apply this result to find $\Lambda^2(V_r)$ in the case where
$p$ is odd, by means of (1.1). Given $\psi^2(V_r)$ it remains only to
calculate $V_r^2$ by the methods available for multiplication in the
Green ring. This
settles a problem left open by Gow and Laffey~\cite{GL} who showed
how to compute $\Lambda^2(V_r)$ when $p=2$.

Our second main result (Theorem 5.1) shows that $\psi^n(V_r)$ has a
strikingly simple form (unlike the much more complicated form that
one gets for $\Lambda^n(V_r)$ or $S^n(V_r)$). Indeed, it turns out
that
$$\psi^n(V_r) = V_{j_1} - V_{j_2} + V_{j_3} - \cdots \pm V_{j_l},$$
where $p^\nu \geq j_1 > j_2 > \cdots > j_l \geq 1$. Thus
the multiplicities of indecomposables in $\psi^n(V_r)$ are only
$0,1$ and $-1$, and the non-zero multiplicities alternate in sign.

The importance of using Adams operations in the study of
$KC$-modules was recognised by Almkvist~\cite{Al}, who studied them
in the case $\nu=1$. An extremely
useful contribution to the study of $\Lambda^n(V_r)$ in the general
case ($\nu \geq 1$) was made by Kouwenhoven~\cite[Theorem 3.5]{Ko},
and his theorem is a key ingredient of our work. By this theorem it
is possible to calculate the values of $\psi_\Lambda^n$ (for all
$n$) on a generating set of $\Rkc$. However, for $n$ not divisible
by $p$, it is known (see \cite{Br}) that $\psi^n$ is an endomorphism
of $\Rkc$.
Thus, in this case, it becomes possible to calculate $\psi^n$ on an
arbitrary element of $\Rkc$. Kouwenhoven studied Adams operations in
his paper~\cite{Ko}, and they also figure in his subsequent
papers~\cite{Ko1, Ko2, Ko3, Ko4}, but his published results seem to
be confined to the case where $\nu=1$.

Hughes and Kemper~\cite{HK} exploited Kouwenhoven's theorem and,
indeed, the results of~\cite[Section 4]{HK} provide, in principle,
a method for
calculating $\Lambda^n(V_r)$ and $S^n(V_r)$ for $n<p$. However, we
believe that our results on Adams operations give a simpler and more
attractive approach.

In a further paper we shall study $\psi_\Lambda^n$ and $\psi_S^n$ on
$\Rkc$ for the general case where $n$ may be divisible by $p$. We
shall prove periodicity results and show that the work of
Symonds~\cite{Sy} may be attractively formulated in terms of Adams
operations.

\vspace{0.4in}
\noindent
{\bf 2. Preliminaries}
\renewcommand{\theequation}{2.\arabic{equation}}
\setcounter{equation}{0}

\noindent Let $G$ be a group and $K$ a field. We consider
$KG$-modules, by which we always mean finite-dimensional right
$KG$-modules, and we write $\Rkg$ for the associated Green ring (or
representation ring). Thus $\Rkg$ is spanned, over $\mathbb{Z}$, by
the isomorphism classes of $KG$-modules and has addition and
multiplication coming from direct sums and tensor products,
respectively. In fact, $\Rkg$ has a $\mathbb{Z}$-basis consisting of
the isomorphism classes of indecomposable $KG$-modules.

For any $KG$-module $V$, we also write $V$ for the corresponding
element of $\Rkg$. Thus, for $KG$-modules $V$ and $W$ we have $V=W$
in $\Rkg$ if and only if $V \cong W$. The elements $V + W$ and $VW$
of $\Rkg$ correspond to $V \oplus W$ and $V \otimes_K W$,
respectively, and the identity element $1$ of $\Rkg$ is the
$1$-dimensional $KG$-module on which $G$ acts trivially. If $V$ is a
$KG$-module and $n$ is a non-negative integer, then we regard
$\Lambda^n(V)$ and $S^n(V)$ as elements of $\Rkg$.

The Adams operations on $\Rkg$ are certain $\mathbb{Z}$-linear maps
from $\Rkg$ to $\Rkg$. We follow the treatment in~\cite{Br}. For
this purpose we need to extend $\Rkg$ to a ring $\mathbb{Q}\Rkg$
where we allow coefficients from $\mathbb{Q}$: thus $\mathbb{Q}\Rkg
\cong \mathbb{Q}\otimes_{\mathbb{Z}} \Rkg$.

For any $KG$-module $V$, define elements of the power-series ring
$\Rkg[[t]]$ by
\begin{eqnarray*}
\Lambda(V,t) &=& 1 + \Lambda^1(V)t +\Lambda^2(V)t^2 + \cdots, \\
S(V,t) &=& 1 + S^1(V)t +S^2(V)t^2 + \cdots .
\end{eqnarray*}
(Since $V$ is assumed to be finite-dimensional, $\Lambda(V,t)$
actually belongs to the polynomial ring $\Rkg[t]$.) Using the formal
expansion of $\log(1+x)$, we have elements $\log\Lambda(V,t)$ and
$\log S(V,t)$ of $\mathbb{Q}\Rkg[[t]]$. Thus we define elements
$\psi_\Lambda^n(V)$ and $\psi_S^n(V)$ of $\mathbb{Q}\Rkg$, for
$n=1,2,\ldots\,$, by the equations
\begin{equation}
\begin{array}{l}
\psi_\Lambda^1(V)t - \frac{1}{2}\psi_\Lambda^2(V)t^2 + \frac{1}{3}
\psi_\Lambda^3(V)t^3 - \cdots = \log\Lambda(V,t),\\[0.15in]
\psi_S^1(V)t + \frac{1}{2}\psi_S^2(V)t^2 + \frac{1}{3}\psi_S^3(V)t^3
+ \cdots = \log S(V,t).
\end{array}
\end{equation}
It is not difficult to prove (for more details see~\cite{Br}) that
$\psi_\Lambda^n(V), \psi_S^n(V) \in \Rkg$ and
$$\psi_\Lambda^n(V+W) = \psi_\Lambda^n(V) + \psi_\Lambda^n(W),\;\;\;
\psi_S^n(V+W) = \psi_S^n(V) + \psi_S^n(W),$$
for all $n\geq 1$ and all $KG$-modules $V$ and $W$. It follows that
the definitions of $\psi_\Lambda^n$ and $\psi_S^n$ may be extended
to give $\mathbb{Z}$-linear functions
$$\psi_\Lambda^n : \Rkg \rightarrow \Rkg,\;\;\; \psi_S^n :
\Rkg \rightarrow \Rkg,$$
called the $n$th Adams operations on $\Rkg$. It is easily verified
that $\psi_\Lambda^1$ and $\psi_S^1$ are equal to the identity map
on $\Rkg$.

For any element $W$ of $\Rkg$ we may now define elements
$\Lambda(W,t)$ and $S(W,t)$ of $\mathbb{Q}\Rkg[[t]]$ by the equations
\begin{eqnarray*}
\Lambda(W,t)&=&\exp(\psi_\Lambda^1(W)t - \textstyle\frac{1}{2}
\psi_\Lambda^2(W)t^2 + \textstyle\frac{1}{3}\psi_\Lambda^3(W)t^3 - \cdots),\\
S(W,t)&=&\exp(\psi_S^1(W)t + \textstyle\frac{1}{2}\psi_S^2(W)t^2 +
\textstyle\frac{1}{3}\psi_S^3(W)t^3 + \cdots).
\end{eqnarray*}
Hence equations (2.1) hold if $V$ is replaced by any element $W$ of
$\Rkg$.

The following result is part of~\cite[Theorem 5.4]{Br}.

\vspace{0.2in} \noindent
{\bf Proposition 2.1.} {\it For every
positive integer\/ $n$ not divisible by the characteristic of\/ $K$,
we have\/ $\psi_\Lambda^n = \psi_S^n$ and each of these maps is a
ring endomorphism of\/ $\Rkg$. Furthermore, under composition of
maps we have}
$$\psi_\Lambda^n \circ \psi_\Lambda^{n'} = \psi_\Lambda^{nn'},
\;\;\;\psi_S^n \circ \psi_S^{n'} = \psi_S^{nn'},$$
{\it for all positive integers\/ $n$ and\/ $n'$ such that\/ $n$ is
not divisible by\/ ${\rm char}\, K$.}

\vspace{0.2in} We shall be mainly concerned with Adams operations
$\psi_\Lambda^n$ and $\psi_S^n$ for $n$ not divisible by
${\rm char}\, K$. For these operations we write $\psi^n$, where
$\psi^n=\psi_\Lambda^n=\psi_S^n$.
We also write $\delta$ for the `dimension' map $\delta : \Rkg
\rightarrow \mathbb{Z}$. This is the $\mathbb{Z}$-linear map
satisfying $\delta(V) = \dim V$ for every $KG$-module $V$.

If $G_1$ is a group of order $1$ then any $KG_1$-module $V$ may be
written as $\delta(V)\cdot 1$ (where $1$ is the identity
element of $R_{KG_1}$) and it is easily verified that
$$\Lambda(V,t)=(1+t)^{\delta(V)}, \;\;\; S(V,t)=(1-t)^{-\delta(V)}.$$
It follows that $\psi_\Lambda^n(V)=\psi_S^n(V)=V$ for all $n$. Thus
each $\psi_\Lambda^n$ and each $\psi_S^n$ is the identity map on
$R_{KG_1}$.

For an arbitrary group $G$ we have homomorphisms $G
\rightarrow G_1$ and $G_1 \rightarrow G$ giving ring homomorphisms
$\alpha : R_{KG_1} \rightarrow \Rkg$ and $\beta : \Rkg \rightarrow
R_{KG_1}$, respectively. Here $\alpha$ is an embedding,
$\beta$ is given by restriction of modules to the identity
subgroup, and $\alpha(\beta(W))=\delta(W)\cdot 1$ for all $W \in
\Rkg$ (where $1$ is the identity element of $\Rkg$). The
formation of exterior and symmetric powers commutes with
restriction: hence $\beta \circ \psi_\Lambda^n =
\psi_\Lambda^n \circ \beta$ and $\beta \circ \psi_S^n =
\psi_S^n \circ \beta$, giving
$$\beta(\psi_\Lambda^n(W))=\beta(\psi_S^n(W))=\beta(W)$$
for all $W\in \Rkg$. On applying $\alpha$ we obtain an equality
of `dimensions':
\begin{equation}
\delta(\psi_\Lambda^n(W))=\delta(\psi_S^n(W))=\delta(W),
\end{equation}
for all $W\in \Rkg$ and all $n \geq 1$.

Now let $p$ be a prime and $K$ a field of characteristic $p$. Let
$\nu$ be a non-negative integer and let $C(p^\nu)$ denote a cyclic
group of order $p^\nu$. It is well known that there are, up to
isomorphism, precisely $p^\nu$ indecomposable $KC(p^\nu)$-modules,
$V_1, V_2, \ldots, V_{p^\nu}$, where $\dim V_r = r$ for $r=1,
\ldots, p^\nu$. (For a proof of this fact see~\cite[Proposition
I.1.1]{AF} or~\cite[Proposition 2.1]{HK}.) Here $V_1$ is the trivial
$1$-dimensional $KC(p^\nu)$-module and $V_{p^\nu}$ is the regular
$KC(p^\nu)$-module.

If $K'$ is an extension field of $K$ there is an embedding
$R_{KC(p^\nu)} \rightarrow R_{K'C(p^\nu)}$ given by extension of
scalars, and the image of $V_r$ is easily seen to be the indecomposable
$K'C(p^\nu)$-module of dimension $r$. Thus $R_{KC(p^\nu)}
\cong R_{K'C(p^\nu)}$. Hence we
regard $R_{KC(p^\nu)}$ as the same for all fields of characteristic
$p$, and write it as $\Rpnu$. The identity element of $\Rpnu$ is
sometimes written as $1$ and sometimes $V_1$.

For each non-negative integer $m$, let $C(p^m)$ be a cyclic group of
order $p^m$ and choose a surjective homomorphism $C(p^{m+1})
\rightarrow C(p^m)$. Thus, for $j \geq m$, the group $C(p^m)$ may be
regarded as a factor group of $C(p^j)$, and there is an injective
homomorphism $\Rpm \rightarrow R_{p^j}$ mapping the $r$-dimensional
indecomposable $KC(p^m)$-module to the $r$-dimensional
indecomposable $KC(p^j)$-module, for $r=1, \ldots, p^m$.

Consequently we may take
$R_{p^0} \subset R_{p^1} \subset \cdots \subset \Rpnu$,
where $\Rpm$ has $\mathbb{Z}$-basis $\{V_1, \ldots, V_{p^m}\}$ for
$m=0, \ldots, \nu$. Throughout the paper we also write $V_0=0$ and
$V_{-r} = -V_r$ for $r = 1,\dots,p^\nu$.

Suppose that $\nu \geq 1$. For $m=0, \ldots, \nu -1$ we define $X_m
\in R_{p^{m+1}}$ by
$$X_m = V_{p^m +1} - V_{p^m - 1},$$
modifying slightly the notation of~\cite{AF}. In particular
$X_0=V_2$. These elements were earlier considered by Green~\cite{Gr}
in a different notation.

\vspace{0.2in} \noindent
{\bf Proposition 2.2.} {\it Let\/ $m \in \{0,
1, \ldots, \nu -1\}$ and\/ $r\in \{0, \ldots, (p-1)p^m\}$. Then}
$$X_m V_r = V_{r+p^m}+V_{r-p^m}.$$

\vspace{0.1in}
\noindent
{\bf Proof.} For $0 < r < (p-1)p^m$ this is given
directly by \cite[(2.3a) and (2.3b)]{Gr}. For $r=0$ it is trivial, and
for $r = (p-1)p^m$ it
follows easily from \cite[(2.3c)]{Gr}. \hfill$\square$

\vspace{0.2in} By the remark immediately after~\cite[Theorem 3]{Gr}
or by~\cite[Proposition I.1.6]{AF}, the Green ring $\Rpnu$ is
generated by the elements $X_0, \ldots, X_{\nu-1}$.

Let $m \in \{0,\dots,\nu\}$. Because $V_{p^m}$ is the
regular $KC(p^m)$-module, we have
$V_{p^m}V_r=rV_{p^m}$ for $r=1, \ldots, p^m$ (by~\cite[VII.7.19
Theorem]{HB}, for example). Hence
\begin{equation}
V_{p^m}W = \delta(W)V_{p^m},
\end{equation}
for all $W \in \Rpm$. It follows that $\mathbb{Z}V_{p^m}$ is an
ideal of $\Rpm$. For $A, B \in \Rpm$ we write
$A \equiv B \;({\rm mod} \; V_{p^m})$
to denote that $A-B\in \mathbb{Z}V_{p^m}$. In fact, such a
congruence gives an equation, by consideration of dimension, namely
$A = B + p^{-m}\delta(A-B)V_{p^m}$.

Note that $V_{p^m}$ is the only projective indecomposable
$KC(p^m)$-module. Also, for $r \in \{1, \ldots, p^m\}$, it is well
known and easy to see that $V_{p^m -r}$ is the Heller translate of
$V_r$ as $KC(p^m)$-module: we write
\begin{equation}
\Omega_{p^m}(V_r) = V_{p^m -r}.
\end{equation}
(For general properties of the Heller translate see~\cite{Be},
for example.)
We extend $\Omega_{p^m}$ to a $\mathbb{Z}$-linear
map $\Omega_{p^m} : \Rpm \rightarrow \Rpm$. Then, for
all $W \in R_{p^m}$, we have
\begin{equation}
\Omega_{p^m}(\Omega_{p^m}(W)) \equiv W \;\;({\rm mod} \; V_{p^m}).
\end{equation}
For $KC(p^m)$-modules $U$ and $V$, consideration of tensor products gives
$$\Omega_{p^m}(UV) \equiv \Omega_{p^m}(U) V \;\;({\rm mod} \; V_{p^m})$$
(see~\cite[Corollary 3.1.6]{Be}). Hence, for all $A,B \in \Rpm$, we have
\begin{equation}
\Omega_{p^m}(AB) \equiv \Omega_{p^m}(A)B \;\;({\rm mod} \; V_{p^m}).
\end{equation}

\vspace{0.3in}
\noindent
{\bf 3. Periodicity and symmetry}
\renewcommand{\theequation}{3.\arabic{equation}}
\setcounter{equation}{0}

\noindent For the remainder of the paper, $p$ is a prime and
$\nu$ is a positive integer. We
consider the Green ring $\Rpnu$ for the cyclic group $C(p^\nu)$
and use the notation of Section 2. In particular,
$X_m = V_{p^m +1} - V_{p^m - 1}$ for $m=0, \ldots, \nu
-1$.

As in~\cite[Section I.1]{AF} and~\cite[Section 4.1]{HK}, let
$\Rpnu$ be extended to a ring
${\widehat R}_{p^\nu}$ generated by $R_{p^\nu}$ and elements
$E_0, \dots, E_{\nu -1}$ satisfying $E_m^2 - X_mE_m + 1 = 0$ for
$m = 0,\dots, \nu -1$. Thus each $E_m$ is invertible in
${\widehat R}_{p^\nu}$ and $X_m = E_m + E_m^{-1}$.
(Note that $E_m$
is written as $\mu_m$ in~\cite{AF} and~\cite{HK}.)

By~\cite[Theorem 3.5]{Ko}, we have
$\Lambda(X_m, t) = 1 + X_m t +t^2$.
Thus
$$
\Lambda(X_m, t) = 1 + (E_m + E_m^{-1})t + t^2 =
(1+E_m t)(1+E_m^{-1}t),
$$
and so, in $(\mathbb{Q}\otimes_{\mathbb{Z}}{\widehat R}_{p^\nu})[[t]]$,
we have
\begin{eqnarray*}
\log\Lambda(X_m,t) &=& \log(1+E_mt) + \log(1+E_m^{-1} t)\\
&=&(E_m + E_m^{-1})t -\textstyle\frac{1}{2}(E_m^2 + E_m ^{-2})t^2 +
\textstyle\frac{1}{3}(E_m^3 + E_m^{-3})t^3 - \cdots.
\end{eqnarray*}
Hence, by (2.1), we obtain
\begin{equation}
\mbox{
$\psi_\Lambda^n(X_m) = E_m^n +E_m^{-n}$ \hspace{0.1in}
for all $n \geq 1$.
}
\end{equation}

For the moment we fix $m$ in the range $0 \leq m \leq \nu-1$ and
write $E=E_m$ and $E^{<n>} = E^n +E^{-n}$ for all
$n \geq 0$. Note that, for $n \geq 1$,
\begin{equation}
E^{<n>}E^{<1>} = E^{<n+1>} + E^{<n-1>}.
\end{equation}
Write $Z=V_{p^m} - V_{p^m-1}$. Thus $Z^2 = 1$, by~\cite[(4.4)]{HK},
and, by \cite[Theorem 4.2]{HK},
$$(ZE-1)((ZE)^{2p-1}-(ZE)^{2p-2} + \cdots + ZE -1) = 0.$$
Since $Z^2 = 1$, we obtain
\begin{equation}
E^{2p}-2ZE^{2p-1}+2E^{2p-2}-\cdots -2ZE +1 =0.
\end{equation}

\vspace{0.2in} \noindent
{\bf Lemma 3.1.} {\it We have\/
$E^{<p+1>}=E^{<p-1>}$.}

\vspace{0.2in} \noindent
{\bf Proof.} Assume first
that $p$ is odd. Multiplying (3.3) by $E^{-p}$ we obtain
$$E^p - 2ZE^{p-1} + \cdots + 2E - 2Z + 2E^{-1} -
\cdots -2ZE^{-(p-1)} + E^{-p} = 0.$$
Hence
\begin{equation}
E^{<p>}=2ZE^{<p-1>}-2E^{<p-2>}+ \cdots -2E^{<1>} +2Z.
\end{equation}
Therefore, by (3.2),
$$E^{<p+1>} +E^{<p-1>} = 2ZE^{<p>} - 2E^{<p-1>}+4ZE^{<p-2>}-\cdots
+ 4ZE^{<1>}-4.$$
Hence, by (3.4), $E^{<p+1>} + E^{<p-1>} = 2E^{<p-1>}$. This gives
$E^{<p+1>}=E^{<p-1>}$, as required. The proof is similar for $p=2$.
\hfill$\square$

\vspace{0.2in} \noindent
{\bf Proposition 3.2.} (i) {\it For\/ $j=0,
\ldots, p$, we have\/ $E^{<2p-j>} = E^{<j>}$.}

(ii) {\it For all\/ $c
\geq 0$, we have\/ $E^{<2p+c>} = E^{<c>}$.}

\vspace{0.2in} \noindent
{\bf Proof.} By Lemma 3.1,
$E^{<p+1>}=E^{<p-1>}$. Multiplying by $E^{<1>}$ we get
$$E^{<p+2>} +E^{<p>} = E^{<p>} + E^{<p-2>}$$
and so $E^{<p+2>}=E^{<p-2>}$. Continuing in this way we obtain
$E^{<p+j>}=E^{<p-j>}$ for $j=0, 1, \ldots, p$. This gives (i).

In particular we have $E^{<2p>}=E^{<0>}=2$. This gives (ii) in the
case $c=0$. Multiplying the equation $E^{<2p>}=2$ by $E^{<1>}$ we
get $E^{<2p+1>} + E^{<2p-1>} =2E^{<1>}$. Since $E^{<2p-1>}=E^{<1>}$,
by (i), we have $E^{<2p+1>}=E^{<1>}$. This gives (ii) in the case
$c=1$. Continuing in this way we get the result for all $c$.
\hfill$\square$

\vspace{0.2in} From now on we write $\psi^n=\psi_\Lambda^n$ for all
$n$ not divisible by $p$. (Thus, in fact,
$\psi^n=\psi_\Lambda^n=\psi_S^n$.)

\vspace{0.2in} \noindent
{\bf Theorem 3.3.} {\it For\/ $j=1, \ldots,
p-1$, we have\/ $\psi^{2p-j} = \psi^j$. Also, if\/ $c$ is any
positive integer not divisible by\/ $p$, we have\/ $\psi^{2p+c} =
\psi^c$.}

\vspace{0.2in} \noindent
{\bf Proof.} As noted in
Section 2, $\Rpnu$ is generated by $\{X_m : 0 \leq m \leq \nu - 1 \}$.
Let $j$ and $c$ be as stated. Then Proposition 3.2 and (3.1) give
$\psi^{2p-j}(X_m)=\psi^j(X_m)$ and
$\psi^{2p+c}(X_m)=\psi^c(X_m)$ for all $m \in \{0,\dots,\nu - 1\}$.
However, by Proposition 2.1, $\psi^{2p-j}$, $\psi^j$, $\psi^{2p+c}$
and $\psi^c$ are endomorphisms of $\Rpnu$. Thus the result follows.
\hfill$\square$

\vspace{0.2in} Let $c$ be any positive integer not divisible by $p$.
Then it is easy to see that there is a unique integer $\gamma(c)$
satisfying the conditions $1 \leq \gamma(c) \leq p-1$ and $c \equiv
\pm \gamma(c) \;\;({\rm mod} \;2p)$. Theorem 3.3 has the following
immediate consequences.

\vspace{0.2in} \noindent
{\bf Corollary 3.4.} {\it For\/ $c$ a
positive integer not divisible by\/ $p$, we have\/
$\psi^c=\psi^{\gamma(c)}$.}

\vspace{0.1in} \noindent
{\bf Corollary 3.5.} {\it Suppose that\/
$p=2$. Then\/ $\psi^c$ is the identity map for every positive
integer\/ $c$ not divisible by\/ $p$.}

\vspace{0.2in}
Let $n$ be a positive integer not divisible by $p$,
and let $m \in \{1, \ldots, \nu\}$. Then
$$V_{p^m-1}^2=(p^m-2)V_{p^m} + V_1,$$
by \cite[(2.5b)]{Gr}.  Hence
\begin{equation}
V_{p^m-1}^n \equiv \left\{
\begin{array}{lll}
V_{p^m-1}&({\rm mod}\;V_{p^m})&\mbox{if $n$ is odd,}\\
V_1&({\rm mod}\;V_{p^m})&\mbox{if $n$ is even.}
\end{array}\right.
\end{equation}

By~\cite[p.\ 362]{Br}, there are $KC(p^\nu)$-modules $Y_d$, for
each divisor $d$ of $n$, such that
\begin{equation}
V_{p^m-1}^n = \sum_{d | n} \phi(d) Y_d,
\end{equation}
where $\phi$ is Euler's function. Also, by~\cite[(4.4) and Theorem
5.4]{Br},
\begin{equation}
\psi^n(V_{p^m-1}) = \sum_{d|n} \mu(d) Y_d,
\end{equation}
where $\mu$ is the M{\"o}bius function.

Note that $\phi(d)=1$ only if $d=1$ or $d=2$. Suppose first that $n$
is odd. Then (3.5) and (3.6) give
$Y_1 \equiv V_{p^m-1}\;({\rm mod}\;V_{p^m})$ and
$Y_d \equiv 0 \;({\rm mod}\;V_{p^m})$ for all $d > 1$.
Thus, by (3.7),
$$\psi^n(V_{p^m-1}) \equiv V_{p^m-1}\;\;({\rm mod}\;V_{p^m}).$$
However, $\delta(\psi^n(V_{p^m-1}))=p^m-1$ by (2.2). Hence
$\psi^n(V_{p^m-1}) = V_{p^m-1}$.

Now suppose that $n$ is even. By (3.5) and (3.6), there
exists $e \in \{1,2\}$ such that
$Y_e \equiv V_1 \;({\rm mod}\; V_{p^m})$ and
$Y_d \equiv 0 \;({\rm mod}\; V_{p^m})$ for all $d \neq e$.
Hence, by (3.7),
$$\psi^n(V_{p^m-1}) \equiv \pm V_1\;\;({\rm mod}\;V_{p^m}).$$
Since $n$ is even, $p \neq 2$.
Thus, using (2.2), we get $\psi^n(V_{p^m-1})=V_{p^m}-V_1$.

Therefore, for all $n$ not divisible by $p$,
\begin{equation}
\psi^n(V_{p^m-1}) = \left\{
\begin{array}{ll}
V_{p^m-1}&\mbox{if $n$ is odd,}\\
V_{p^m}-V_1&\mbox{if $n$ is even.}
\end{array}\right.
\end{equation}
By similar, but much easier, arguments we obtain
\begin{equation}
\psi^n(V_{p^m})=V_{p^m}\;\;\;\mbox{for all $n$ not divisible by $p$.}
\end{equation}

By \cite[(2.5b)]{Gr}, we have
\begin{equation}
V_{p^m-1}V_r = (r-1)V_{p^m} + V_{p^m-r},
\end{equation}
for $r=1,\dots,p^m$. (Recall that $V_0 = 0$.)
Hence, by Proposition 2.1 and (3.9),
\begin{equation}
\psi^n(V_{p^m-1})\psi^n(V_r)=(r-1)V_{p^m} + \psi^n(V_{p^m-r}),
\end{equation}
for all $n$ not divisible by $p$. Note that (3.8)--(3.11) hold,
trivially, for $m=0$. Thus they hold for all $m \in \{0,\dots,\nu\}$.

\vspace{0.2in} \noindent
{\bf Proposition 3.6.} {\it Let\/ $n$ be an even positive integer
not divisible by\/ $p$} ({\it thus\/ $p$ is odd\/}),
{\it and let\/ $m \in \{0, \ldots, \nu\}$. Then, for\/ $r=1, \ldots,
p^m$, we have}
$$\psi^n(V_r) + \psi^n(V_{p^m-r}) = V_{p^m}.$$

\vspace{0.2in} \noindent
{\bf Proof.} By (3.8) and (3.11),
$$(V_{p^m}-V_1)\psi^n(V_r) = (r-1)V_{p^m} + \psi^n(V_{p^m-r}).$$
However, $V_{p^m}\psi^n(V_r)=rV_{p^m}$ by (2.2) and (2.3). This
gives the required result. \hfill$\square$

\vspace{0.2in} By (3.10) and (2.4) we have, for all $W \in R_{p^m}$,
\begin{equation}
V_{p^m-1}W \equiv \Omega_{p^m}(W)\;\;({\rm mod} \; V_{p^m}).
\end{equation}

\vspace{0.2in} \noindent
{\samepage
{\bf Proposition 3.7.} {\it Let\/ $n$ be an
odd positive integer not divisible by\/ $p$, and let\/ $m \in \{0,
\ldots, \nu\}$. Then, for\/ $r=1, \ldots, p^m$, we have}
$$\psi^n(V_{p^m-r})\equiv \Omega_{p^m}(\psi^n(V_r))
\;\;({\rm mod}\;V_{p^m}).$$

\vspace{0.2in} \noindent
{\bf Proof.} By (3.8),
$\psi^n(V_{p^m-1}) = V_{p^m-1}$. Hence, by (3.11),
$$V_{p^m-1}\psi^n(V_r) = (r-1)V_{p^m} + \psi^n(V_{p^m-r}).$$
Thus the result follows by (3.12). \hfill$\square$
}

\vspace{0.2in}
Propositions 3.6 and 3.7 are partial generalisations
of~\cite[Propositions 5.4(d) and 5.4(e)]{Al}. Stronger results will
be given below in Corollary 5.2.

We conclude this section by showing that, when $n=2$, Proposition
3.6 implies Gow and Laffey's ``reciprocity theorem'' \cite[Theorem
1]{GL}. This may be stated in the Green ring as follows (after
correction of the obvious misprint in~\cite{GL}).

\vspace{0.2in} \noindent
{\bf Corollary 3.8.} {\it Let\/
$p$ be odd and\/ $m \in \{1, \ldots, \nu\}$. Then, for\/ $r=1,
\ldots, p^m$,}

(i) $\;\Lambda^2(V_r) = (r-\frac{1}{2}(p^m+1))V_{p^m} +
S^2(V_{p^m-r})$,

(ii) $S^2(V_r) = (r-\frac{1}{2}(p^m-1))V_{p^m} +
\Lambda^2(V_{p^m-r})$.

\vspace{0.2in} \noindent
{\bf Proof.} Since (i) and (ii) are essentially the same we
prove only (i). It is well known
that $S^2(V_{p^m-r}) + \Lambda^2(V_{p^m-r}) = V_{p^m-r}^2$.
Thus
$$\Lambda^2(V_r) - S^2(V_{p^m-r}) = \Lambda^2(V_r) +
\Lambda^2(V_{p^m-r}) - V_{p^m-r}^2.$$
By (1.1) (which follows from (2.1)), we have
$\Lambda^2(V_r) = \textstyle\frac{1}{2}(V_r^2 - \psi^2(V_r))$;
and a similar statement holds for $\Lambda^2(V_{p^m-r})$. Hence
$$\Lambda^2(V_r) - S^2(V_{p^m-r}) = \textstyle\frac{1}{2}(V_r^2 -
V_{p^m-r}^2) - \textstyle\frac{1}{2}(\psi^2(V_r) +
\psi^2(V_{p^m-r})).$$
However, by (2.5), (2.6) and (2.4), we have
$$V_r^2 \equiv \Omega_{p^m}(\Omega_{p^m}(V_r^2)) \equiv
(\Omega_{p^m}(V_r))^2 \equiv V_{p^m-r}^2 \;\;({\rm mod} \; V_{p^m}),$$
so that
$V_r^2 - V_{p^m-r}^2 = (2r-p^m)V_{p^m}$.
Also, we have $\psi^2(V_r) + \psi^2(V_{p^m-r}) = V_{p^m}$, by
Proposition 3.6. Thus
$$\Lambda^2(V_r) - S^2(V_{p^m-r}) = \textstyle\frac{1}{2}
(2r-p^m)V_{p^m} - \textstyle\frac{1}{2}V_{p^m},$$
which gives the required result. \hfill$\square$

\vspace{0.4in} \noindent
{\bf 4. Recursion}
\renewcommand{\theequation}{4.\arabic{equation}}
\setcounter{equation}{0}

{\samepage
\noindent
Define elements $g_0(t), g_1(t), \ldots$ of
$\mathbb{Z}[t]$ by $g_0(t)=2$, $g_1(t) = t$ and, for $n \geq 2$,
\begin{equation}
g_n(t) = tg_{n-1}(t) - g_{n-2}(t).
\end{equation}
The $g_n(t)$ can be seen to be Dickson
polynomials of the first kind, and can be given by an explicit
formula, but we do not need this.
}

\vspace{0.2in} \noindent
{\bf Proposition 4.1.} {\it For\/ $n \geq
1$ and\/ $m \in \{0, \ldots, \nu-1\}$, we have}
$$\psi_\Lambda^n(X_m) = g_n(X_m).$$

\vspace{0.2in} \noindent
{\bf Proof.} Clearly
$\psi_\Lambda^1(X_m) = X_m$ and, by (3.1),
$\psi_\Lambda^2(X_m) = X_m^2 - 2$. Hence the result holds for
$n \leq 2$.
It is easy to check from (3.1) and (3.2) that, for $n\geq 3$,
$$\psi_\Lambda^n(X_m) = X_m\psi_\Lambda^{n-1}(X_m) -
\psi_\Lambda^{n-2}(X_m).$$
Thus the result follows by induction and (4.1). \hfill$\square$

\vspace{0.2in}
For $n < p$, Proposition 4.1 can be deduced from (3.1) and
\cite[(I.1.4) and (I.1.5)]{AF}. Our next result is a
reformulation of \cite[Lemma 4.2]{F}, but we give a proof for the
convenience of the reader.

\vspace{0.2in} \noindent
{\bf Proposition 4.2.} {\it Let\/ $m \in
\{0, \ldots, \nu-1\}$,  $r \in \{1, \ldots, p^m\}$, and\/ $i \in
\{0, \ldots, p-1\}$. Then}
$$g_i(X_m)V_r = V_{ip^m+r}-V_{ip^m-r}.$$

\vspace{0.2in} \noindent
{\bf Proof.} The result is
clear for $i=0$ because, by convention, $V_{-r}$ denotes $-V_r$.
Since $g_1(X_m)=X_m$, the result for $i=1$ is given by Proposition
2.2. Now suppose that $2 \leq i \leq p-1$ and the result holds
for $i-1$ and $i-2$. Then, by (4.1) and the inductive hypothesis,
\begin{eqnarray*}
g_i(X_m)V_r &=& X_m g_{i-1}(X_m)V_r - g_{i-2}(X_m)V_r\\
&=&X_m(V_{(i-1)p^m+r} - V_{(i-1)p^m-r}) - (V_{(i-2)p^m+r} -
V_{(i-2)p^m-r}).
\end{eqnarray*}
It is easy to verify that $(i-1)p^m+r$ and $(i-1)p^m-r$ belong to
$\{0, \ldots, (p-1)p^m\}$. Hence, by Proposition 2.2,
{\samepage
\begin{eqnarray*}
g_i(X_m)V_r &=& (V_{ip^m+r} + V_{(i-2)p^m+r}) - (V_{ip^m-r} +
V_{(i-2)p^m-r})\\
&&\;\;\; \mbox{} - (V_{(i-2)p^m+r} - V_{(i-2)p^m-r})\\
&=&V_{ip^m+r} - V_{ip^m-r},\end{eqnarray*} as required.
\hfill$\square$
}

\vspace{0.2in}
For a positive integer $c$ not divisible by $p$, let
$\gamma(c)$ be as defined in Section~3. Note that $1 \leq
\gamma(c) \leq p-1$.

\vspace{0.2in} \noindent
{\bf Corollary 4.3.} {\it Let\/ $m \in \{0,
\ldots, \nu-1\}$. For\/ $r \in \{1, \dots, p^m\}$ and\/
$c$ any positive integer not divisible by\/ $p$, we have}
$$\psi^c(X_m)V_r = V_{\gamma(c)p^m+r} - V_{\gamma(c)p^m-r}.$$

\vspace{0.2in} \noindent
{\bf Proof.} By Corollary
3.4, $\psi^c(X_m) = \psi^{\gamma(c)}(X_m)$. Hence, by Proposition
4.1, $\psi^c(X_m)=g_{\gamma(c)}(X_m)$. Thus the result follows by
Proposition 4.2. \hfill$\square$

\vspace{0.2in} For $m \in \{0, \ldots, \nu-1\}$ and $i \in \{1,
\ldots, p-1\}$ let $\theta_{ip^m} : \Rpm \rightarrow R_{p^{m+1}}$ be the
$\mathbb{Z}$-linear map defined by
\begin{equation}
\theta_{ip^m}(V_r) = V_{ip^m+r} - V_{ip^m-r}
\end{equation}
for $r=1, \ldots, p^m$. Corollary 4.3 gives the
following result.

\vspace{0.2in} \noindent
{\bf Corollary 4.4.} {\it Let\/ $m \in \{0,
\ldots, \nu-1\}$. Let\/ $c$ be any positive integer not divisible
by\/ $p$ and let\/ $W \in \Rpm$. Then
$$\psi^c(X_m)W = \theta_{\gamma(c)p^m}(W).$$}

\vspace{0.01in}
Define elements $f_{-1}(t), f_0(t), f_1(t), \ldots$
of $\mathbb{Z}[t]$ by $f_{-1}(t)=0, f_0(t)=1, f_1(t)=t$ and, for $n
\geq 2$,
\begin{equation}
f_n(t) = tf_{n-1}(t) - f_{n-2}(t).
\end{equation}
The $f_n(t)$ can be seen to be Dickson
polynomials of the second kind, and can be given an explicit
formula, but we do not need this. The following result is
straightforward to prove by induction.

{\samepage
\vspace{0.2in} \noindent
{\bf Lemma 4.5.} {\it For all\/ $n \geq
0$,}
$$f_n = \left\{
\begin{array}{ll}
g_n + g_{n-2} + \cdots + g_3 +g_1 &\mbox{if $n$ is odd,}\\
g_n + g_{n-2} + \cdots + g_2 + 1&\mbox{if $n$ is even.}
\end{array}\right.$$

\vspace{0.2in} Our next result is essentially the same
as~\cite[Lemma 6]{McF}, but we
give a proof for the convenience of the reader.
}

\vspace{0.2in} \noindent
{\bf Proposition 4.6.} {\it Let\/ $m \in
\{0, \ldots, \nu-1\}$. Then, for\/ $r \in \{1, \ldots, p^m\}$
and\/ $k \in \{0, \ldots, p-1\}$, we have
$$V_{kp^m+r} = f_k(X_m)V_r + f_{k-1}(X_m)V_{p^m-r}.$$}

\vspace{0.01in} \noindent
{\bf Proof.} We use
induction on $k$. The result is clear for $k=0$. It is true for
$k=1$ because
$V_{p^m+r} = X_mV_r + V_{p^m-r}$ by Proposition 2.2.

Now suppose that $k \in \{2, \ldots, p-1\}$ and that the result is
true for $k-1$ and $k-2$. By (4.3), the inductive hypothesis, and
Proposition 2.2, we obtain
\begin{eqnarray*}
f_k(X_m)V_r + f_{k-1}(X_m)V_{p^m-r}
&=& X_m(f_{k-1}(X_m)V_r + f_{k-2}(X_m)V_{p^m-r})\\
&& \;\;\; \mbox{} - (f_{k-2}(X_m)V_r + f_{k-3}(X_m)V_{p^m-r})\\
&=& X_m V_{(k-1)p^m +r} - V_{(k-2)p^m+r}\\
&=& V_{kp^m+r},
\end{eqnarray*}
as required. \hfill$\square$

\vspace{0.2in}
In the statement of the main result of this section
it is convenient to extend the definition of $\gamma$ by setting
$\gamma(0)=0$. Recalling that $\theta_{ip^m}$ is defined by (4.2)
for $i \in \{1,\dots,p-1\}$,
we also define $\theta_0$ to be the identity map
on $R_{p^m}$.

\vspace{0.2in} \noindent
{\bf Theorem 4.7.} {\it Let\/ $m \in  \{0,
\ldots, \nu-1\}$ and let\/ $n$ be a positive integer not divisible
by\/ $p$. Let\/ $s$ be a positive integer satisfying\/ $p^m < s \leq
p^{m+1}$ and write\/ $s=kp^m+r$, where\/ $1 \leq r \leq p^m$ and\/
$1 \leq k \leq p-1$. Then
$$\psi^n(V_s) = \sum_{j \in \{0,\dots,k\}
\atop j \equiv k \,({\rm mod}\, 2)}
\theta_{\gamma(jn)p^m}(\psi^n(V_r)) \; +
\sum_{j \in \{0,\dots,k\}
\atop j \not\equiv k \,({\rm mod}\, 2)}
\theta_{\gamma(jn)p^m}(\psi^n(V_{p^m-r})).$$}

\vspace{0.2in} \noindent
{\bf Proof.} By Proposition 4.6, we have
$V_s = f_k(X_m)V_r + f_{k-1}(X_m)V_{p^m-r}$.
Suppose first that $k$ is odd. Then, by Lemma 4.5 and Proposition
4.1, we obtain
\begin{eqnarray*}
V_s&=& (\psi^k + \psi^{k-2} + \cdots + \psi^1)(X_m) V_r\\
&& \;\;\; \mbox{} + (\psi^{k-1} + \psi^{k-3} + \cdots + \psi^2)(X_m)
V_{p^m-r} + V_{p^m-r}.
\end{eqnarray*}
By Proposition 2.1 it follows that
\begin{eqnarray*}
\psi^n(V_s)&=& (\psi^{kn} + \psi^{(k-2)n} + \cdots + \psi^n)(X_m)
\psi^n(V_r)\\
&& \;\;\; \mbox{} + (\psi^{(k-1)n} + \psi^{(k-3)n} + \cdots +
\psi^{2n})(X_m) \psi^n(V_{p^m-r}) + \psi^n(V_{p^m-r}).
\end{eqnarray*}
Therefore, by Corollary 4.4,
\begin{eqnarray*}
\psi^n(V_s)&=& (\theta_{\gamma(kn)p^m} + \theta_{\gamma((k-2)n)p^m}
+ \cdots + \theta_{\gamma(n)p^m})(\psi^n(V_r))\\
&& \;\;\; \mbox{} + (\theta_{\gamma((k-1)n)p^m} + \theta_{\gamma((k-3)n)p^m}
+ \cdots + \theta_{\gamma(2n)p^m} + \theta_0)(\psi^n(V_{p^m-r})),
\end{eqnarray*}
as required. The proof for even $k$ is similar. \hfill$\square$

\vspace{0.2in}
Theorem 4.7 allows us to calculate $\psi^n(V_s)$
for all $s$, and for all $n$ not divisible by $p$, by
elementary arithmetic and without
the need for multiplication in $R_{p^\nu}$.

For example, take $p=7$ and $\nu=2$. Let us
calculate $\psi^4(V_{23})$. Thus $n=4$ and $s=23$.
In order to apply Theorem 4.7 we take $m=1$ and write
$23=3\cdot 7+2$. (Thus $k=3$ and $r=2$.) It is easy to check
that $\gamma(4)=4$, $\gamma(2\cdot4)=6$ and $\gamma(3\cdot4)=2$.
Thus, by Theorem 4.7,
\begin{equation}
\psi^4(V_{23}) = (\theta_{28} + \theta_{14})(\psi^4(V_2)) +
(\theta_{0} + \theta_{42})(\psi^4(V_5)).
\end{equation}
We next calculate $\psi^4(V_2)$, writing $2= 1 \cdot 1 + 1$ in
order to use Theorem 4.7. Thus
$$\psi^4(V_2) = \theta_4(\psi^4(V_1)) + \theta_{0}(\psi^4(V_0))
= \theta_{4}(V_1)
= V_5 - V_3.$$
We can calculate $\psi^4(V_5)$ in a similar way, or by means of
Proposition 3.6, to obtain
$\psi^4(V_5)= V_7 - V_5 + V_3$. Thus, by (4.4),
{\samepage
\begin{eqnarray*}
\psi^4(V_{23}) &=& (\theta_{28} + \theta_{14})(V_5 - V_3) +
(\theta_{0} + \theta_{42})(V_7 - V_5 + V_3)\\
&=& (V_{33} - V_{23}) + (V_{19} - V_9) - (V_{31} - V_{25}) -
(V_{17} - V_{11})\\
&& \;\;\; \mbox{} +  V_7 + (V_{49} - V_{35}) - V_5 - (V_{47} - V_{37}) +
V_3 + (V_{45} - V_{39})\\
&=& V_{49} - V_{47} + V_{45} - V_{39} + V_{37} - V_{35} + V_{33} -
V_{31} + V_{25}\\
&& \;\;\; \mbox{} - V_{23} + V_{19} - V_{17} + V_{11}
- V_{9} + V_7 - V_{5} + V_3.
\end{eqnarray*}
}
We see that the indecomposables occurring have all subscripts of the
same parity and have multiplicities that alternate between $+1$ and
$-1$, in decreasing order of subscript. It turns out that these
statements hold in general. We shall prove them in Theorem 5.1 in
the next section.

\vspace{0.4in} \noindent
{\bf 5. The form of $\psi^n(V_s)$}
\renewcommand{\theequation}{5.\arabic{equation}}
\setcounter{equation}{0}

\vspace{0.2in} \noindent
{\bf Theorem 5.1.} {\it Let\/ $n$ be a
positive integer not divisible by\/ $p$, and let\/ $s \in \{1,
\ldots, p^\nu\}$. Write\/ $\lambda(s)$ for the smallest non-negative
integer such that\/ $s \leq p^{\lambda(s)}$.}

(i) {\it There are integers\/ $j_1,\dots,j_l$ such that\/
$p^{\lambda(s)} \geq j_1 > j_2 > \cdots > j_l \geq 1$ and
$$\psi^n(V_s) = V_{j_1} - V_{j_2} + V_{j_3}- \cdots \pm V_{j_l}.$$}

\vspace{-0.2in}
(ii) {\it If \/ $n$ is even\/} ({\it so that\/ $p$ is odd\/})
{\it then\/ $j_1, \ldots, j_l$ are odd.
If\/ $n$ is odd then\/ $j_1, \ldots, j_l$ have the same parity as\/
$s$.}

\vspace{0.2in}
Before giving the proof we derive an improvement of
Propositions 3.6 and 3.7.

\vspace{0.2in} \noindent
{\bf Corollary 5.2.} {\it Let\/ $n$ be a
positive integer not divisible by\/ $p$, and let\/ $s \in \{1,
\ldots, p^m\}$, where\/ $m \in \{0,\dots,\nu\}$.}

(i) {\it If\/ $n$ is even then one of\/ $\psi^n(V_s)$ and\/
$\psi^n(V_{p^m-s})$ has the form
$$V_{j_1} - V_{j_2} + \cdots \pm V_{j_l}$$
and the other has the form
$$V_{p^m} - V_{j_1} + V_{j_2} - \cdots \mp V_{j_l},$$
where\/ $j_1,\dots,j_l$ are odd and\/
$p^m > j_1 > j_2 > \cdots > j_l \geq 1$.}

\vspace{0.1in}
(ii) {\it If \/ $n$ is odd then\/ $\psi^n(V_s)$ and\/
$\psi^n(V_{p^m-s})$ have the forms
$$\psi^n(V_s)= V_{j_1} - V_{j_2} +  V_{j_3} - \cdots + V_{j_l},$$
$$\psi^n(V_{p^m-s})= V_{p^m-j_l} -  \cdots + V_{p^m-j_3} - V_{p^m-j_2}
+ V_{p^m-j_1},$$
where\/ $l$ is odd, $j_1,\dots,j_l$ have the parity of\/ $s$,
and\/ $p^m \geq j_1 > j_2 > \cdots > j_l \geq 0$.}

\vspace{0.2in} \noindent
{\bf Proof.} (i) This is
immediate from Theorem 5.1 and Proposition 3.6.

(ii) If $p=2$ then $\psi^n$ is the identity map, by Corollary 3.5,
and the result is clear. Thus we may assume that $p$ is odd. We
argue according to the parity of $s$.

Suppose first that $s$ is odd. By Theorem 5.1 we may write
$$\psi^n(V_s) = V_{j_1} - V_{j_2} + V_{j_3} - \cdots \pm V_{j_l},$$
where $j_1, \ldots, j_l$ are odd and $p^m \geq j_1 > j_2 > \cdots >
j_l \geq 1$. By (2.2),
$$\delta(V_{j_1} - V_{j_2} + \cdots \pm V_{j_l}) =s.$$
Since $s$ is odd it follows that $l$ must be odd, and so $\psi^n(V_s)$
has the required form.
By Theorem 5.1, $\psi^n(V_{p^m-s})$ is a linear combination of terms
$V_i$ where $i$ has the parity of $p^m-s$; so
$\psi^n(V_{p^m-s})$ does not involve $V_{p^m}$. Thus,
by Proposition 3.7,
$$\psi^n(V_{p^m-s}) = V_{p^m-j_l} - \cdots + V_{p^m-j_3}-V_{p^m-j_2}+
V_{p^m-j_1}.$$
(Note here that we may have $p^m-j_1=0$.) Thus the result holds for
$s$ odd. If $s$ is even then $p^m-s$ is odd and we may interchange
the roles of $V_s$ and $V_{p^m-s}$ in the above argument.
\hfill$\square$

\vspace{0.2in} \noindent
{\bf Proof of Theorem 5.1.}
For each integer $a$ let $[a]$ denote the congruence class of $a$
modulo $2$ and let $R[a]$ denote the additive subgroup of
$R_{p^\nu}$ spanned by all $V_i$ with $[i] = [a]$. Thus
$R[a] = R[0]$ or $R[a] = R[1]$. Observe
that (i) and (ii) of Theorem 5.1 are equivalent to (i) and
the statement that $\psi^n(V_s) \in R[ns+n+1]$.

To prove the theorem we use induction on $m$, where $m = \lambda(s)$.
Since $\psi^n(V_1) = V_1$, statements (i) and (ii) are trivial for
$m=0$. Let $m < \nu$ and assume that (i) and (ii) hold for all $s$
with $\lambda(s) \leq m$. Now take $s$ such that $\lambda(s) = m+1$.
We shall prove that (i) and (ii) hold for $V_s$.
Write $q = p^m$, so that $pq = p^{m+1}$. Also, write $s = kq+r$,
where $1 \leq r \leq q$ and $1 \leq k \leq p-1$, as in Theorem 4.7.
Thus $\psi^n(V_r)$ and $\psi^n(V_{q-r})$ are covered by
the inductive hypothesis.

For each non-negative integer $a$ define $U_a$ by
$$U_a= \left\{
\begin{array}{ll}
\psi^n(V_r)\; & \mbox{if $\;[a] = [k]$,}\\
\psi^n(V_{q-r})\; & \mbox{if $\;[a] \neq [k]$.}
\end{array}\right.$$
Then, by Theorem 4.7,
\begin{equation}
\psi^n(V_s) = \sum_{j=0}^k \theta_{\gamma(jn)q} (U_j).
\end{equation}

We have
$\psi^n(V_r) \in R[nr+n+1]$ and
$\psi^n(V_{q-r}) \in R[n(q-r)+n+1]$,
by the inductive hypothesis. It follows easily that
$U_j \in R[nr+n+1+(j+k)nq]$
for $j=0,\dots,k$. By the definition of
$\theta_{\gamma(jn)q}$ (see (4.2)), we obtain
$$\theta_{\gamma(jn)q}(U_j) \in R[nr+n+1+(j+k)nq+\gamma(jn)q].$$
However, $[\gamma(jn)] = [jn]$. Thus
$$\theta_{\gamma(jn)q}(U_j) \in
R[n(kq+r)+n+1] = R[ns+n+1].$$
Hence, by (5.1), we have $\psi^n(V_s) \in R[ns+n+1]$. Thus it remains
only to prove that (i) holds. We deal separately with the cases
where $n$ is even and $n$ is odd.

Suppose first that $n$ is even, so that $p$ is odd. Clearly
$\lambda(pq-s) \leq m+1$. Also, by
Proposition 3.6,
$\psi^n(V_s) + \psi^n(V_{pq-s}) = V_{pq}$.
It follows that if (i) holds for $V_{pq-s}$ then it
holds for $V_s$. Thus, by the inductive hypothesis, we may assume
that $\lambda(pq-s) = m+1$.
Either $s < \frac{1}{2}pq$ or $pq-s <\frac{1}{2}pq$.
Therefore, without loss of generality, we may assume
that $s < \frac{1}{2}pq$.

Since $s=kq+r < \frac{1}{2}pq$, we have $k \leq
\frac{1}{2}(p-1)$.
Suppose that $\gamma(s_1n)=\gamma(s_2n)$,
where $s_1, s_2 \in \{1, \ldots, k\}$. Then
$s_1n \equiv \pm s_2n \;({\rm mod}\;2p)$.
Since $p \nmid n$ we obtain
$s_1 \equiv \pm s_2\; ({\rm mod}\;p)$.
Hence $s_1 \mp s_2 \equiv 0 \; ({\rm mod}\;p)$.
However,
$s_1, s_2 \in \{1, \ldots, \frac{1}{2}(p-1)\}$
because $k \leq \frac{1}{2}(p-1)$.
Therefore $s_1=s_2$.
Thus the numbers $\gamma(n), \gamma(2n), \ldots,
\gamma(kn)$ are distinct. They are even, since $n$ is even. Hence
we may write
$$\{\gamma(n), \gamma(2n), \ldots, \gamma(kn)\}
= \{a_1, a_2, \ldots, a_k\},$$
where the $a_j$ are even and
$p-1 \geq a_1 > a_2 > \cdots > a_k \geq 2$.
Also, set $a_{k+1} = 0$.

By (5.1) we have
\begin{equation}
\psi^n(V_s) = \theta_{a_1q}(W_1) + \cdots + \theta_{a_kq}(W_k) +
\theta_{a_{k+1}q}(W_{k+1}),
\end{equation}
where $W_j \in \{\psi^n(V_r), \psi^n(V_{q-r})\}$ for each $j$.

For integers $a$ and $b$ with $pq \geq a \geq b \geq 0$, let
$M[a,b]$ denote the set of all elements $Y$ of $\Rpq$ that can be
written in the form
$$Y=V_{i_1}-V_{i_2} +V_{i_3} - \cdots + V_{i_{h-1}} - V_{i_h}$$
where $h$ is even and $a \geq i_1 \geq i_2 \geq \cdots \geq i_h \geq
b$.
To prove (i) it suffices to show that $\psi^n(V_s) \in M[pq,0]$, for
then we obtain the required expression for $\psi^n(V_s)$ by
cancellation and by removal of terms $V_0$.

Suppose that $pq \geq c_1 \geq c_2 \geq \cdots \geq c_{d+1} \geq 0$
and $Y_j \in M[c_j,c_{j+1}]$ for $j=1, \ldots, d$. Then, clearly,
$Y_1 + Y_2 + \cdots + Y_d \in M[c_1, c_{d+1}]$.

By the inductive hypothesis, each $W_j$ belongs to $M[q,0]$, since
we may introduce a term $V_0$ if necessary to give even length to
the expression for $W_j$. It follows easily that
$\theta_{a_jq}(W_j)$ belongs to $M[(a_j+1)q, (a_j-1)q]$, for $j=1,
\ldots, k$. Hence
$$\theta_{a_jq}(W_j) \in M[(a_j+1)q,(a_{j+1}+1)q],$$
for $j=1,\ldots, k$, because $a_j \geq a_{j+1} + 2$. Also,
$$\theta_{a_{k+1}q}(W_{k+1}) =
W_{k+1} \in M[q,0] = M[(a_{k+1}+1)q,0].$$
Therefore, by (5.2), we have
$\psi^n(V_s) \in M[(a_1+1)q,0] \subseteq M[pq,0]$,
as required.

We now turn to the remaining case, and assume that
$n$ is odd.

Since Theorem 5.1 holds for $V_r$ and $V_{q-r}$,
by the inductive
hypothesis, Corollary 5.2(ii) holds for $V_r$ and $V_{q-r}$. Thus we
may write
\begin{equation}
\psi^n(V_r) = V_{j_1} - V_{j_2} + V_{j_3} - \cdots + V_{j_l},
\end{equation}
\begin{equation}
\psi^n(V_{q-r}) = V_{q-j_l}- \cdots + V_{q-j_3} - V_{q-j_2} +
V_{q-j_1},
\end{equation}
where $l$ is odd and
$q \geq j_1 > j_2 > \cdots > j_l \geq 0$.

Suppose that $\gamma(s_1n)=\gamma(s_2n)$,
where $s_1, s_2 \in \{1, \ldots, k\}$. Then we have
$s_1n \equiv \pm s_2n \;({\rm mod}\;2p)$.
Since $n$ is coprime to $2p$ we obtain
$s_1 \equiv \pm s_2\;({\rm mod}\;2p)$.
Hence $s_1 \mp s_2 \equiv 0 \; ({\rm mod}\;2p)$. Since
$s_1,s_2 \in \{1,\dots,p-1\}$, it follows that $s_1=s_2$.
Consequently, the numbers $\gamma(n), \gamma(2n), \ldots,
\gamma(kn)$ are distinct and we may write
$$\{\gamma(n), \gamma(2n), \ldots, \gamma(kn)\} =
\{a_1, a_2, \ldots, a_k\},$$
where $p-1 \geq a_1 > a_2 > \cdots > a_k \geq 1$.
Also, set $a_{k+1} = 0$.

Since $n$ is odd, we have
$[\gamma(jn)] = [j]$, and (5.1) may be written
\begin{equation}
\psi^n(V_s) = \theta_{a_1q}(U_{a_1}) + \theta_{a_2q}(U_{a_2}) +
\cdots + \theta_{a_kq}(U_{a_k}) + \theta_{a_{k+1}q}(U_{a_{k+1}}).
\end{equation}

With $j_1,\dots,j_l$ as in (5.3) and (5.4), define $T_{aq}$,
for each $a \in \{0,\dots,p-1\}$, by
\vspace{0.1in}
$$T_{aq} = \left\{
\begin{array}{ll}
V_{aq+j_1} - V_{aq+j_2} + \cdots + V_{aq +j_l}\; & \mbox{if
$\;[a] = [k]$,}\\
V_{aq+q-j_l} - \cdots -V_{aq+q-j_2} + V_{aq + q-j_1}\; & \mbox{if
$\;[a] \neq [k]$.}
\end{array}\right.$$

\vspace{0.05in}
\noindent
Then it can be checked that
$\theta_{aq}(U_a) = T_{aq} - T_{(a-1)q}$
for all $a \in \{1,\dots,p-1\}$. Also, $\theta_{0q}(U_0) = T_{0q}$.
Thus, by (5.5),
$$\psi^n(V_s) = T_{a_1q} - T_{(a_1-1)q} + T_{a_2q} - T_{(a_2-1)q}
+ \cdots + T_{a_kq} - T_{(a_k-1)q} + T_{a_{k+1}q}.$$
If $a_j-1=a_{j+1}$ for some $j \in \{1, \ldots, k\}$ then we may
cancel two adjacent terms in this expression.
After all such cancellations we obtain
\begin{equation}
\psi^n(V_s) = T_{b_1q} - T_{b_2q} +
\cdots - T_{b_{d-1}q} + T_{b_dq},
\end{equation}
where $d$ is odd and $p-1 \geq b_1 > b_2 > \cdots > b_d \geq 0$.

For integers $a$ and $b$ where $pq \geq a \geq b \geq 0$, let
$N[a,b]$ denote the set of all elements $Y$ of $\Rpq$ that can be
written in the form
$$Y = V_{i_1}-V_{i_2} +\cdots - V_{i_{h-1}} + V_{i_h}$$
where $h$ is odd and $a \geq i_1 \geq i_2 \geq \cdots \geq i_h \geq
b$. To prove (i) it suffices to show that $\psi^n(V_s) \in N[pq,0]$.

By the definition of $T_{aq}$ we see that
$T_{b_jq} \in N[(b_j+1)q,b_jq]$
for $j=1,\ldots, d$. However, $b_jq \geq (b_{j+1}+1)q$, for $j \leq
d-1$. Therefore, by (5.6),
$$\psi^n(V_s) \in N[(b_1+1)q,b_dq] \subseteq N[pq,0],$$
as required. \hfill$\square$

\end{document}